\documentclass[12pt,reqno]{amsart}
\usepackage{amsmath, amsfonts, amssymb, amsthm, hyperref}
\usepackage{bm}
\allowdisplaybreaks[4]
\def\l{\left}
\def\r{\right}
\def\bg{\bigg}
\def\({\bg(}
\def\){\bg)}
\def\t{\text}
\def\f{\frac}

\def\ls{\leqslant}

\def\bi{\binom}

\def\eq{\equiv}

\def\Z{\mathbb Z}

\def\N{\mathbb N}

\def\<{\langle}
\def\>{\rangle}
\def\1{{\bf 1}}

\theoremstyle{plain}
\newtheorem{theorem}{Theorem}

\newtheorem{lemma}{Lemma}
\newtheorem{corollary}{Corollary}
\theoremstyle{definition}
\newtheorem*{Ack}{Acknowledgment}
\theoremstyle{remark}

\newtheorem{Data Availability}{Data Availability}

\makeatletter

\newcommand{\Rmnum}[1]{\expandafter\@slowromancap\romannumeral #1@}
\makeatother
\numberwithin{equation}{section}

\begin{document}
\title[Legendre symbols related to $D_p(b,1)$]{Legendre symbols related to $D_p(b,1)$}
\author
[] {Xin-Qi Luo and Wei Xia}

\address {(Xin-Qi Luo) Department of Mathematics, Nanjing
University, Nanjing 210093, People's Republic of China}
\email{lxq15995710087@163.com}

\address{(Wei Xia) Department of Mathematics, Nanjing
University, Nanjing 210093, People's Republic of China}
\email{wxia@smail.nju.edu.cn}

\begin{abstract}
Let $p$ be an odd prime. For any $b,c\in\mathbb{Z}$, Z.-W. Sun introduced the new-type determinant
$$D_p(b,c)=|(i^2+bij+cj^2)^{p-2}|_{1\leqslant i,j\leqslant p-1},$$
and studied its arithmetic properties. In this paper we mainly prove that
$$\left(\frac{D_p(b,1)}{p}\right)=\left(\frac{2b}{p}\right)$$
when $(\frac{b^2-4}{p})=-1$ and $p\equiv1\pmod 4$. As an application of our result, we confirm several conjectures of Sun.
\end{abstract}
\keywords{determinants, Legendre symbol, Lucas sequence.}
\subjclass[2020]{Primary 11C20; Secondary 11A07, 11B39.}
\thanks{Supported by the National Natural Science Foundation of China (Grant No. 12371004).}
\maketitle

\section{Introduction}
\setcounter{lemma}{0} \setcounter{theorem}{0}
\setcounter{equation}{0}\setcounter{proposition}{0}
\setcounter{Rem}{0}\setcounter{conjecture}{0}
For an $n\times n$ matrix $M=[a_{ij}]_{1\leqslant i,j\leqslant n}$ over a commutative ring, we denote $|M|$ as the determinant of $M$.

Let $p$ be an odd prime, and let $(\frac{.}{p})$ be the Legendre symbol.
Z.-W.Sun \cite{S19} studied some determinants with the form $(\f{i^2+cij+dj^2}p)$,
where $c,d\in\Z$. For example, Sun showed that for any prime $p\eq3\pmod4$, the $p$-adic number
$$2\left|\frac{1}{i^2+j^2}\right|_{1\leqslant i,j \leqslant (p-1)/2}
$$
is a quadratic residue modulo $p$.

In 2002, Wu, She, and Ni \cite{WSN} proved Z.-W.Sun's \cite{S19} conjecture: For each prime $p\equiv 2\pmod{3}$
 $$2\left|\frac{1}{i^2-ij+j^2}\right|_{1\leqslant i,j \leqslant p-1}$$
 is a quadratic residue modulo $p$.

Let $p$ be an odd prime. For any $c,d\in\mathbb{Z}$, Z.-W. Sun \cite{S22} introduced the determinant
$$D_p(b,c):=|(i^2+bij+cj^2)^{p-2}|_{1\leqslant i,j\leqslant p-1},$$
and studied the Legendre symbol $(\frac{D_p(b,c)}{p})$. When $i^2+bij+cj^2\not\eq0\pmod p$, it can be obtained from Fermat's little theorem that
$$
(i^2+bij+cj^2)^{p-2}\eq\f1{i^2+bij+cj^2}\pmod p.
$$
Sun \cite[(1.7)]{S22} pointed out that
$$D_p(-b,c)\eq\l(\f{-1}p\r)D_p(b,c)\pmod p.$$
Based on the Wu-She-Ni result, it shows that
$$D_p(1,1)=\l(\f{(-1)^{(p-1)/2}D_p(-1,1)}p\r)=\l(\f{-2}p\r)$$
for $p\eq2\pmod3$.

Given two integers $A$ and $B$, the Lucas sequence $u_n=u_n(A,B)$ $(n\in\mathbb{N})$ and its companion $v_n=v_n(A,B)$ $(n\in\mathbb{N})$ are defined as follows:
$$u_0=0,\ u_1=1,\ \mbox{and}\ u_{n+1}=Au_n-Bu_{n-1}\ (n=1,2,3,\ldots),$$
and
$$v_0=2,\ v_1=A,\ \mbox{and}\ v_{n+1}=Av_n-Bv_{n-1}\ (n=1,2,3,\ldots).$$

Recently, Luo and Sun\cite{LS22} investigated $(\frac{D_p(b,c)}{p})$ systematically. In fact, they use trinomial coefficient and Lucas sequence deduced the following result.

\begin{lemma}\label{Lem-U}
Let $p$ be an odd prime, and let $b,c\in\Z$ with $p\nmid c(b^2-4c)$. Let
\begin{equation}\label{U}
U(k)= \bi{p-2}k_{b,c}+c^{p-1-k}\bi{p-2}{p-1-k}_{b,c}.
\end{equation}

{\rm (i)} If $U(k)\equiv 0\pmod{p}$ for some $k\in\{2,\ldots,p-2\}$, then
$$
\l(\f{D_p(b,c)}{p}\right)=0.
$$

{\rm (ii)}
If $U(k)\not\equiv 0\pmod{p}$ for all $2\ls k\ls p-2$, then
\begin{equation}\label{D-p}\begin{aligned}
&\left(\f cp\r)^{(p-1)(p-3)/8}\l(\f{D_p(b,c)}{p}\right)\\
=&\left(\frac{4c-b^2+2c(\f{b^2-4c}p)}p\r)
\l(\f{2cu_{p-1}(-b,c)-b}p\r)\l(\f{U(p-2)U(\frac{p-1}{2})}{p}\right).\\
\end{aligned}\end{equation}
\end{lemma}

For an odd prime $p$, we let $\mathbb Z_p$ denote the ring of all $p$-adic integers. It is well known that
$$
\mathbb Q\cap \mathbb Z_p=\left\{\frac ab: a,b\in \mathbb Z\ \t{and}\ p\nmid b\right\}.
$$
For any $x\in \mathbb Z_p$, there is an unique $r\in \{0,\ldots ,p-1\}$ with $x\equiv r\pmod p$, and we define $(\frac xp)$ as $(\frac rp)$.

Based on Lemma 1.1, Y.-F. She and H.-L. Wu\cite{SW22} obtained the following result.
They pointed out that the most important case is the case $c=1$ and the other case is related to it. Specifically, they obtained the following result.

\begin{theorem}[She-Wu\cite{SW22}]Let $p$ be an odd prime. For any $b\in\mathbb{Z}_p$, the following results hold.\\
{\rm (i)} Suppose $(\frac{b^2-4}{p})=1$.
If $v_k(-b,1)\equiv0\pmod p$ for some $k=0,\ldots ,p-3$, then $\l(\frac{D_p(b,1)}{p}\r)=0$, otherwise
$$\l(\frac{D_p(b,1)}{p}\r)=\l(\frac{-b-2}{p}\r)^{(p-1)/2}.$$
{\rm (ii)} Suppose $(\frac{b^2-4}{p})=-1$. If $v_k(-b,1)\equiv0\pmod p$ for some $k=2,\ldots ,p-1$, then $\l(\frac{D_p(b,1)}{p}\r)=0$, otherwise $$\l(\frac{D_p(b,1)}{p}\r)=\l(\frac{-b-2}{p}\r)^{(p-1)/2}\l(\frac{2b}{p}\r).$$
{\rm (iii)} For any $c\in \mathbb{Z}_p$, we have
$$\l(\frac{D_p(b,c^2)}{p}\r)=\l(\frac{c}{p}\r)^{(p-1)/2}\l(\frac{D_p(b/c,1)}{p}\r).$$
\end{theorem}

Motivated by a lemma of Sun \cite[Lemma 2.2]{S12} and the method of J. C. Lagarias \cite{L85}, we obtain the following further result.

\begin{theorem}Let $p$ be a prime with $p\equiv 1\pmod{4}$, and let $b\in \mathbb Z$ with $(\frac{b^2-4}{p})=-1$. Then we have
$$\left(\frac{D_p(b,1)}{p}\right)=\left(\frac{2b}{p}\right).$$
\end{theorem}

\begin{corollary}Let $p$ be a prime with $p\equiv 1\pmod{4}$. For any $b\in \mathbb Z$ with $(\frac{b^2+4}{p})=-1$, we have
$$\left(\frac{D_p(b,-1)}{p}\right)=\left(\frac{b}{p}\right).$$
\end{corollary}

We also confirm Conjecture 4.5 in \cite{S22}.
\begin{theorem}Let $p$ be an odd prime. If $(\frac{p}{7})=-1$ and $p\not\equiv 15\pmod{16}$, we have
$$\left(\frac{D_p(1,16)}{p}\right)=\left(\frac{-2}{p}\right).$$
\end{theorem}

We need the following lemma \cite[Lemma 2.2]{S12} on quadratic algebraic integers.

\begin{lemma}[Z.-W. Sun \cite{S12}]
Let $A,B\in \mathbb Z$ and $\Delta=A^2-4B$.
For any odd prime $p\nmid B\Delta$, we have
$$
\l(\frac {A\pm\sqrt {\Delta}}2\r)^{p-(\frac \Delta p)}\equiv B^{(1-(\frac \Delta p))/2}\pmod p
$$
in the ring of algebraic integers.
\end{lemma}

\section{Proofs of Theorems}
\setcounter{lemma}{0} \setcounter{theorem}{0}
\medskip
Our proof of Theorem 1.2 is motivated by Lemma 1.2 and the method of Lagarias \cite{L85}.
\medskip

\noindent{\it Proof of Theorem 1.2}.
We define $v_n=v_n(-b,1)\ (n\in\N)$ as follows:
$$v_0=2,\ v_1=-b,\ \t{and}\ v_{n+1}=-bv_{n}-v_{n-1}\ \t{for}\ n=1,2,3,\ldots.$$The characteristic equation $x^2+bx+1=0$ has two roots
$$\alpha=\frac{-b+\sqrt{b^2-4}}{2}\ \mbox{and}\ \beta=\frac{-b-\sqrt{b^2-4}}{2}.$$
It view of Binet's formula, for any $k\in\mathbb{N}$ we have
$$v_k=\alpha^k+\beta^k.$$
Note that $\alpha \beta=1$ and hence
$$\frac{\alpha}{\beta}=\frac{\alpha^2}{\alpha \beta}=\alpha^2.$$
Observe that
\begin{equation}\label{equivalent}
v_k\equiv0\pmod p\Leftrightarrow(\frac{\alpha}{\beta})^k\equiv-1\pmod p\Leftrightarrow(\alpha^2)^k\equiv-1\pmod p,
\end{equation}
where the congruences are in the ring $O_K$ of algebraic integers in $K=\mathbb{Q}[\sqrt{b^2-4}]$. Since $(\frac{b^2-4}{p})=-1$, we know that $(p)$ is a prime ideal in $O_K$ and $O_K/(p)$ is a finite field of $p^2$ elements. For any nonzero $\gamma\in O_K$, we use $\mathrm {ord}_p\gamma$ to denote the least positive integer $n$ such that $\gamma^n\equiv1\pmod p$.

According to Theorem 1.1 (ii), it suffices to show that $v_k\not\equiv0\pmod p$ for $2\leqslant k\leqslant p-1$ provided that $p\equiv 1\pmod{4}$ and $(\frac{b^2-4}{p})=-1$.

Let $\theta=\alpha^2$.
Then
$$
\theta^{\frac {p+1} 2}=\alpha^{p+1}\equiv 1\pmod p
$$
by Lemma 1.2, and hence $\mathrm {ord}_p{\theta}\mid(p+1)/2$. Since $p\equiv1\pmod 4$, $\mathrm {ord}_p{\theta}$ must be odd.


 Suppose that $\theta^{n}\equiv-1\pmod p$ for a positive integer $n$. Then $\theta^{2n}\equiv1\pmod p$, which implies that ord$_p{\theta}\mid 2n$. Since ord$_p{\theta}$ is odd, we obtain  ord$_p{\theta}\mid n$, i.e., $\theta^{n}\equiv1\pmod p$. This leads a contradiction.

By the above, we claim that $\alpha^{2k}= \theta^k\not\equiv-1\pmod p$ for any positive integer $k$. In view of \eqref{equivalent}, we get $v_k\not\equiv0\pmod p$ for any positive integer $k$.
So far we have completed the proof of Theorem 1.2. \qed

\medskip

\noindent{\it Proof of Corollary 1.1}.
Let $q=\frac {p-1}2!$. As $p\equiv1\pmod 4$, by Wilson's theorem we have
$$
\l(\frac {p-1}2!\r)^2\equiv \prod_{k=1}^{\frac{p-1}2}k(p-k)=(p-1)!\equiv -1\pmod p.
$$
 In view of Theorem 1.1 (iii), we have
$$\left(\frac{D_p(b,-1)}{p}\right)=\left(\frac{D_p(b/q,1)}{p}\right).$$
It is easy to check that $$\left(\frac{(b/q)^2-4}{p}\right)=\left(\frac{-b^2-4}{p}\right)=\left(\frac{b^2+4}{p}\right)=-1.$$
Thus, with the help of Theorem 1.2, we obtain
$$\left(\frac{D_p(b/q,1)}{p}\right)=\left(\frac{2b/q}{p}\right)=\left(\frac{2bq}{p}\right).$$
By \cite[Lemma 2.3]{S19}, we have
$$\left(\frac{2q}{p}\right)=1.$$
So we have the desired equality
$$
\l(\frac{D_p(b,-1)}p\r)=\l(\frac bp\r).
$$
\qed

The following three corollaries were originally conjectured by Z.-W. Sun\cite[Conjecture 4.2-4.4]{S22}.

\begin{corollary}For any prime $p\equiv 1\pmod 4$ with $p\equiv \pm2 \pmod 5$, we have
$$\left(\frac{D_p(1,-1)}{p}\right)=1.$$
\end{corollary}
This follows form Corollary 1.1 with $b=1$.
\begin{corollary}For any odd prime $p$, we have
$$\left(\frac{D_p(2,-1)}{p}\right)=-1\Leftrightarrow p\equiv 5\pmod 8.$$
\end{corollary}
\proof
By Sun \cite[Theorem 1.3]{S22}, we have $\left(\frac{D_p(2,-1)}{p}\right)=0$ for any prime $p\equiv 3\pmod 4$.
When $p\equiv 1\pmod 8$, by Theorem 1.1, we have $\left(\frac{D_p(2,-1)}{p}\right)=1\ \t{or}\ 0$. By Corollary 1.1 we have
$$\left(\frac{D_p(2,-1)}{p}\right)=\l(\frac 2p\r)=-1$$
when $p\equiv 5\pmod 8$. This concludes the proof.
\begin{corollary}For any prime $p\equiv \pm2 \pmod 5$, we have
$$\l(\f{D_p(3,1)}p\r)=\begin{cases}(\f 6p)&\t{if}\ p\eq1\pmod4,
\\0&\t{if}\ p\eq3\pmod4.\end{cases}$$
\end{corollary}
\proof
When $p\equiv 1\pmod 4$, we have $(\f 5p )= \l(\f p5 \r)=-1$, then by Theorem 1.2 we have
$$
\left(\frac{D_p(3,1)}{p}\right)=\l(\f 6p\r).
$$
When $p\equiv 3\pmod 4$ and $p\equiv \pm2 \pmod 5$, by \cite[Corollary 3.1]{S03} we have $v_{\f {p+1}2}(-3,1)\equiv0\pmod p$,
thus $p\mid D_p(3,1)$ by Theorem 1.1 (ii). This concludes the proof.\qed
\medskip

\noindent{\it Proof of Theorem 1.3.}
According to Theorem 1.1 (3), we have
 $$\l(\frac{D_p(1,16)}{p}\r)=\l(\frac{D_p(4^{-1},1)}{p}\r).$$

 {\it Case 1}. $p\equiv 1\pmod{4}$.

 Note that
 $$\l(\frac{4^{-2}-4}{p}\r)=\l(\frac {-7}p\r)=\l(\frac 7p\r)=\l(\frac p7\r)=-1.$$
Applying Corollary 1.1 we get
$$\l(\frac{D_p(4^{-1},1)}{p}\r)=\l(\frac{2}{p}\r)=\l(\frac{-2}{p}\r).$$

 {\it Case 2}. $p\eq 3\pmod 4$, but $p\not\eq 15\pmod {16}$.

 According to Theorem 1.1 (ii), if $v_k(-4^{-1},1)\not\equiv0\pmod p$, for all $2\leq k<p-1,$ then $$\l(\frac{D_p(4^{-1},1)}{p}\r)=\l(\frac{-4^{-1}-2}{p}\r)^{(p-1)/2}\l(\frac{2^{-1}}{p}\r)=\l(\frac{-2}{p}\r).$$
It is easy to see that $(-4)^kv_k(-4^{-1},1)=v_k(1,16)=v_{2k}(-3,4)$ from Binet's formula. So it suffice, to prove that $v_{2k}(-3,4)\not\equiv0\pmod p$ for any positive integer $h$.

The characteristic equation $x^2+3 x+4=0$ has two roots
$$\alpha=\f {-3+\sqrt {-7}}{2}\ \t{and}\ \beta=\f {-3-\sqrt {-7}}{2}.$$
It view of Binet's formula, for any $k\in\mathbb{N}$ we have
$$v_{2k}(-3,4)=\alpha^{2k}+\beta^{2k}\in \mathbb Z\l[\frac{-1+\sqrt {-7}}{2}\r].$$
As $(\frac {-7}p)=-1,$ by Lemma 1.3 we have
$$
\alpha^{p+1}\equiv 4\pmod p
$$
in the ring $\mathbb Z[\frac{-1+\sqrt {-7}}{2}].$

Let $\theta =(\frac {\alpha}2)^2=\frac {\alpha}{\beta}$. Then
$$
\theta^{\frac {p+1}2}=\frac {\alpha^{p+1}}{2^{p+1}}\equiv 1\pmod p,
$$
and hence
$$\mathrm {ord}_p{\theta}\mid \frac {p+1}2.$$
Note that
$$
v_{2k}(-3,4)\equiv 0\pmod p\Leftrightarrow(\frac{\alpha}{\beta})^{2k}\equiv -1\pmod p\Leftrightarrow {\theta}^{2k}\equiv -1\pmod p.
$$
If ${\theta}^{2k}\equiv -1\pmod p$, then $\mathrm {ord}_p{\theta}\mid 4k$ and hence $\mathrm {ord}_p{\theta}$ divides $(\frac {p+1}2,4k)$.
If $p\equiv 3\pmod 8$, then $(\frac {p+1}2,4k)\mid 2k$ and hence ${\theta}^{2k}\not\equiv -1\pmod p.$

Now we consider the case $p\equiv 7\pmod {16}$. As $\alpha=(\frac {1+\sqrt{-7}}2)^2,$ we have
\begin{align*}
{\theta}^{\frac {p+1}4}&=\l(\frac {\alpha}2\r)^{\frac {p+1}2}\equiv \frac{( \frac 2p)}2\l(\frac{1+\sqrt{-7}}2\r)^{p+1}\equiv \frac {1+\sqrt{-7}}{8}(1+\sqrt{-7})^p\\
&=\frac {1+(\sqrt{-7})^{p+1}+\sqrt{-7}(1+(\sqrt{-7})^{p-1})}8\equiv \frac 18\l(1+(-7)^{\frac {p+1}2}+\sqrt{-7}\l(1+\l(\frac {-7}p\r)\r)\r)\\
&\equiv \frac 18\l(1-7\l(\frac {-7}p\r)\r)=1\pmod p.
\end{align*}
and hence $\mathrm {ord}_p{\theta}$ divides $\frac {p+1}4\equiv 2\pmod 4$. As $(\frac {p+1}4,4k)\mid 2k$, we cannot have ${\theta}^{2k}\equiv -1\pmod p.$
In view of the above, we have completed the proof.
\begin{Ack}
The authors are grateful to Z.-W. Sun for his valuable suggestions.
\end{Ack}
\noindent{\bf Data Availability Statement.} Data sharing is not applicable to this article as no new data were created or analyzed in this study.
\\

\noindent{\bf Compliance with Ethical Standards.} This study was funded by the National Natural Science Foundation of China (Grant No. 12371004). All authors disclosed no relevant relationships. This article does not contain any studies with human participants or animals performed by any of the authors.

\end{document}